\documentclass[11pt]{amsart}
\usepackage{comment}
   
\usepackage[usenames]{color}
\usepackage{amssymb}
\usepackage{latexsym}
\usepackage{graphics} 
\usepackage{graphicx}
\usepackage{epstopdf}  
\usepackage{enumitem}
\usepackage[all]{xy}
\usepackage{verbatim} 
\DeclareSymbolFont{AMSb}{U}{msb}{m}{n}
\DeclareMathSymbol{\N}{\mathbin}{AMSb}{"4E}
\DeclareMathSymbol{\Z}{\mathbin}{AMSb}{"5A}
\DeclareMathSymbol{\R}{\mathbin}{AMSb}{"52}
\DeclareMathSymbol{\Q}{\mathbin}{AMSb}{"51}
\DeclareMathSymbol{\I}{\mathbin}{AMSb}{"49}
 
\newtheorem{theorem}{Theorem}

\newtheorem{cor}{Corollary}

\newtheorem{lemma}{Lemma}

\newtheorem{defi}{Definition}
\newtheorem{q}[theorem]{Question}

\newtheorem*{menasco}{Menasco's Theorem}
\newtheorem*{menascoState}{Menasco's Statement}

\newtheorem*{wrongM}{Sawollek's Theorem}

\newtheorem*{AdamsM}{Adams et al.'s Theorem}

\newtheorem*{futer}{Futer, Kalfagianni, and Purcel's Theorem}

\newtheorem*{ozawa}{Ozawa's Theorem}
\newtheorem*{This}{Thistlethwaite's Theorem}
\newtheorem*{cromwell}{Cromwell's Theorem}
\newtheorem*{nontriv}{Kauffman, Murasugi, and Thistlethwaite's Theorem}
\newtheorem*{Stoimenow}{Stoimenow's Theorem}

\theoremstyle{definition}

\newcommand{\bi}{\begin{itemize}}
\newcommand{\ei}{\end{itemize}}
\newcommand{\be}{\begin{enumerate}}
\newcommand{\ee}{\end{enumerate}}

\begin{document}

\begin{abstract} Menasco proved the surprising result that if $G$ is a reduced, alternating, connected projection of a link $L$ and $G$ is prime then $L$ is prime.  This result has been generalized to other classes of links, tangles, and spatial graphs.  We draw attention to some issues with previous splitting results about tangles and spatial graphs, and obtain new more general results for tangles and spatial graphs.
\end{abstract}
\title{Splittings of Tangles and Spatial Graphs}
  \date{\today}
  
  \author[E. Flapan, H.N.\ Howards]{Erica Flapan, Hugh Howards}

    
    \address{Emerita Professor, Department of Mathematics, Pomona College, Claremont, CA 91711, USA}

\address{Department of Mathematics, Wake Forest University, Winston-Salem, NC 27109, USA}

\thanks{The first author was supported in part by NSF Grant DMS-1607744}

  \maketitle
  
\section {Introduction}

One would not expect to be able to conclude much about a link by looking at a single projection.  Yet in 1984, Menasco proved the striking result that if a link has a reduced alternating projection, then you can determine from such a projection whether the link is split or composite \cite{Menasco}.  We are interested in analogous results for tangles and spatial graphs.   Menasco's result and the relevant definitions are given below.

\begin{menasco}  Let $L$ be a link with a reduced alternating projection $G$.  Then
\begin{enumerate}
\item If $G$ is connected, then $L$ is non-split.
\item If $L$ is non-split and $G$ is prime, then $L$ is prime.
\end{enumerate}
\end{menasco}

\begin{defi}  Let $L$ be a link with projection $G$ on a sphere $P$.
	\begin{itemize}
	
	\item We say $G$ is {\bf  connected} if the shadow obtained by replacing the crossings of $G$ with singularities is connected. 
		\item We say $L$ is {\bf prime} if every sphere meeting $L$ transversely in two points bounds a ball meeting $L$ in an unknotted arc.  
		\item We say $G$ is {\bf prime} if every circle in $P$ meeting $G$ transversely in two points bounds a disk meeting $G$ in an arc with no crossings.  
		
	\end{itemize}
\end{defi}

Part (1) and (2) of Menasco's Theorem have been generalized to larger classes of links by Tsukamoto \cite{Tsukamoto}, Cromwell \cite{CromHomo}, Thistlethwaite \cite{ThistlethwaiteAdequate}, Ozawa \cite{Ozawa3} and Futer, Kalfagianni, Purcell \cite{Futer}, as well as to tangles by Cromwell~\cite{Crom}, and to spatial graphs by Sawollek~\cite{Sawollek2} and Adams et al.~\cite{Adams}.  Unfortunately, the study of projections for tangles and spatial graphs has not always been precise.  The goals of the current paper are to draw attention to some of the issues in past results, and to further generalize Menasco's Theorem for tangles and spatial graphs in Theorems~\ref{tangle}, \ref{alternating}, and \ref{figureeight}.  

In particular, in Section~\ref{Links}, we state some generalizations to broader classes of links, and mention issues with the statement of Menasco's Theorem.  Then, in Section~\ref{tangles}, we discuss Cromwell's generalization of Menasco's result to tangles and use results for links to prove a generalization to a class of tangles larger than those considered by Cromwell.   In Section~\ref{Graphs}, we discuss definitions and results of Sawollek and Adams et al. for spatial graphs.   Finally, in Section~\ref{OurGraphs}, we build on our result for tangles to obtain results for spatial figure eight graphs that go beyond those of Sawollek and Adams et al..

\section{Generalizations of Menasco's Theorem to other links}\label{Links}

 The statement of Part (2) of Menasco's Theorem that we gave in the introduction is not, in fact, how he stated his result.  Rather, he stated his theorem as follows.

\begin{menascoState} \cite{Menasco} \label{MenWrong} Let $L$ be a link with a reduced alternating projection $G$.  Then
\begin{enumerate}
\item If $G$ is connected, then $L$ is non-split.
\item Suppose that $L$ is non-split.  Then $L$ is prime if for every disk $D$ in the plane of projection whose boundary meets $G$ transversely in just two points, $G\cap D$ is a trivial arc.  
\end{enumerate}
\end{menascoState}

However, as stated, the hypothesis in Part (2) is never satisfied for a non-trivial projection $G$ of a link.  In particular, let $D$ be a disk that includes all of $G$ except for a small arc with no crossings.  Then $G\cap D$ will necessarily be non-trivial.  Thus, Part (2) of the above theorem can only be applied to a link with no crossings.  The same mistake occurs in some other statements and definitions related to generalizations of Menasco's Theorem.  Thus, we take Menasco's Theorem to be as we stated it in the Introduction.  Here we present some related results for link projections.

The following theorem was proved independently by Kauffman, Murasugi, and Thistlethwaite in 1987.

\begin{nontriv} \cite{Kauff}, \cite{Mura}, \cite{Thistle}   Any link with a projection which has at least one crossing and is connected, reduced, and alternating, is non-trivial.
\end{nontriv}



  In 2002, Ozawa generalized Menasco's Theorem to \emph{positive} links, that is those which have an oriented projection in which all crossings are positive. 
  
 \begin{ozawa} \cite{Ozawa3} Let $L$ be a link whose projection $G$ is connected and positive.  Then $L$ is non-split.  Furthermore, if $L$ is non-trivial and $G$ is prime, then $L$ is prime.
 \end{ozawa}

 In 2003, Stoimenow proved the following generalization of Kauffman, Murasugi, and Thistlethwaite's Theorem, which together with Ozawa's result shows that any link which has a reduced positive connected projection is non-trivial.

 \begin{Stoimenow}\cite{Stoi}  Any knot which has a reduced positive projection with at least one crossing is non-trivial.  
 \end{Stoimenow}

Part (1) of Menasco's Theorem was generalized to semi-adequate links by Thistlethwaite~\cite{This} in 1988, and Part (2) was generalized to semi-adequate links  meeting an additional condition by Futer et al.~\cite{Futer} in 2015.  Recall the following definitions.

\begin{figure}[h!]
 	\centering	\includegraphics[width=.45\textwidth]{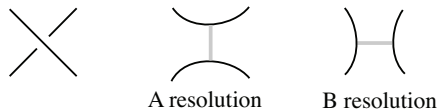}
 	\caption{If we add these grey segments to the all-$A$ resolution or all-$B$ resolution we get a $3$-valent graph.}
 	\label{smoothings}
 \end{figure}

   \begin{defi}  Let $G$ be a projection of a link.  If we resolve all of the crossings as in the picture in the center of Figure~\ref{smoothings} (ignoring the grey segments), we obtain a collection of circles which we call the {\bf all-$A$ resolution} of $G$.  If we resolve all of the crossings as in the picture on the right, we obtain a collection of circles which we call the {\bf all-$B$ resolution} of $G$.  \end{defi}

 \begin{defi} Consider the $3$-valent planar graph obtained by adding grey edges to an all-$A$ or all-$B$ resolution of a link projection $G$ as in Figure~\ref{smoothings}.  We say $G$ is {\bf $A$-adequate} or {\bf $B$-adequate} if every such grey segment has its endpoints on distinct circles of the all-$A$ resolution or all-$B$ resolution, respectively.  In either case, we say that $G$ is {\bf semi-adequate}.  If $G$ is both $A$-adequate and $B$-adequate, then we say $G$ is {\bf adequate}.
 \end{defi}			
 
\begin{This}\cite{This} Let $L$ be a link whose projection $G$ is semi-adequate and has at least one crossing.
\begin{enumerate}
\item  If $G$ is reduced, then $L$ is non-trivial.
\item  If $G$ is connected, then $L$ is non-split.
\end{enumerate}
\end{This}

To state the result of Futer, Kalfagianni, and Purcell, we need some additional definitions.

 
 A {\bf twist region} of a link projection is a maximal row of bigons in the projection.  For example, all of the crossings in the standard projection of a $(2,n)$-torus link are in a single twist region, while a standard projection of a $(p,q,r)$-pretzel link with $|p|$, $|q|$, $|r|>1$ has three twist regions.  

\begin{defi}  We say an $A$-adequate (resp. $B$-adequate) link projection is {\bf twist $A$-adequate} (resp. {\bf twist $B$-adequate}) if it has the property that any pair of grey arcs going between the same pair of circles in an all-$A$ (resp. all-$B$) resolution come from resolving crossings in the same twist region.  We say a link projection is {\bf twist-adequate} if it is either twist $A$-adequate or twist $B$-adequate.   
 \end{defi}
 
 \begin{futer}\cite{Futer}  Let $L$ be a link whose projection $G$ is connected, reduced, and twist-adequate.  Then $L$ is prime if and only if $G$ is prime. \end{futer}
 
Any alternating link projection is semi-adequate  \cite{Menasco} and any positive link projection is semi-adequate \cite{Ozawa3}.   So if it weren't for the condition on twist regions, then the results of Futer et al. would imply both Menasco's Theorem and Ozawa's Theorem.

\section{Splittings of tangles}\label{tangles}

Menasco's Theorem has also been generalized to $2$-string tangles $(B,T)$ whose projections $(\beta,\tau)$ are {\bf strongly alternating} (i.e., such that $\tau$ has at least one crossing and its closures $N(\tau)$ and $D(\tau)$ are reduced and alternating).  Note that some authors require a strongly alternating tangle projection to be connected, but Cromwell \cite{Crom} does not.  
 
 
 \begin{defi}\label{tangleDef} Let $(B,T)$ be a non-rational $2$-string tangle.   
 	\begin{itemize}
 		\item We say $(B,T)$ is {\bf prime} if every sphere in $B$ meeting $T$ transversely in two points bounds a ball in $B$ whose intersection with $T$ is an unknotted arc.  Otherwise, we say $(B,T)$ is {\bf composite}.  
 		
 		\item  We say the projection $(\beta,\tau)$ of $(B,T)$ is {\bf prime} if every circle in $\beta$ meeting $\tau$ transversely in two points bounds a disk in $\beta$ whose intersection with $\tau$ is an arc with no crossings.  Otherwise, we say $(\beta,\tau)$ is {\bf composite}.
 	\end{itemize}
 \end{defi}
 
 Observe that if $(B,T)$ is split (i.e., there is a disk which separates the strings), then either it is rational or there is a sphere in $B$ meeting $T$ in two points bounding a ball in $B$ whose intersection with $T$ is a knotted arc.  Hence, every split $2$-string tangle is either rational or composite.   
  
 \begin{cromwell}  {(Theorem 5(b) of \cite{Crom})}  Let $(B,T)$ be a tangle whose projection $(\beta,\tau)$ is strongly alternating. Then $(B, T)$ is prime if and only if $(\beta,\tau)$ is prime.
 \end{cromwell}
 
\begin{figure}[h!]
 	\centering	\includegraphics[width=.5\textwidth]{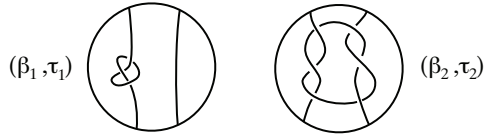}
 	\caption{$(\beta_1,\tau_1)$ is composite but its numerator closure is prime. $(\beta_2,\tau_2)$ is strongly alternating but is a sum of tangles which are not strongly alternating.}
 	\label{NotPrime}
 \end{figure}

 Cromwell's proof contained some gaps which we illustrate in Figure~\ref{NotPrime}.  In particular,
his proof asserts that if one of the closures of a tangle is a prime knot, then the tangle must also be prime.  However,  $(\beta_1,\tau_1)$ is composite while its numerator closure is a prime knot.  Cromwell's proof also asserts that if a strongly alternating tangle projection can be written as a tangle sum, then each of the summands is strongly alternating.  However, the strongly alternating projection $(\beta_2,\tau_2)$ can be split into a sum of two tangles which are not strongly alternating.   Though it isn't necessary for Cromwell's proof, it is also worth noting that $(\beta_2,\tau_2)$ is a prime tangle while its closure $D(\tau)$ is a composite knot.  Thus the primeness or compositeness of a tangle is not determined by the primeness or compositeness of its closures.
 
 Shimokawa \cite{Shimo} asserts the same theorem as Cromwell under the additional hypothesis that $(\beta,\tau)$ is connected, saying simply that the result can be proved using Menasco's machinery.

\medskip
 
 Our result for tangles is stronger than Cromwell's Theorem, and implies it as a corollary.  We begin with two lemmas.


 \begin{lemma}
 	Let $(\beta, \tau)$ be a connected projection of a $2$-string tangle such that $D(\tau)$ is reduced, and let $(\beta',\tau')$ to be a $\pi$-rotation of $(\beta, \tau)$.  Then $N(\tau+ \tau')$ is reduced and connected, and if $D(\tau)$ is alternating or positive, or $N(\tau)$ and $D(\tau)$ are twist-adequate, then $N(\tau+ \tau')$ is alternating, positive, or twist-adequate, respectively.

 	\label{lem:cp}
 \end{lemma}

 \begin{proof}  Let $P$ be the sphere of projection, let $\gamma=\partial \beta$, $\gamma'=\partial\beta'$, and let $\beta_1$ and $\beta_2$  be the strips joining $(\beta,\tau)$ and $(\beta',\tau')$ illustrated in Figure~\ref{diskpicture}. 
 	
 	\begin{figure}[h!]
 		\centering	\includegraphics[width=.5\textwidth]{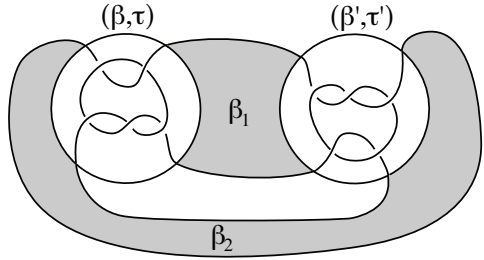}
 		\caption{ $N(\tau  + \tau' )$ with the strips $\beta_1$ and $\beta_2$ indicated.}
 		\label{diskpicture}
 	\end{figure}

 \medskip

 	\noindent{\bf Step 1:  We show that $N(\tau  + \tau' )$ is reduced and connected.}

 	Suppose for the sake of contradiction that there exists a circle $C\subseteq P$ which is either a separating circle for $N(\tau  + \tau' )$ or meets $N(\tau  + \tau' )$ transversely in a single crossing.  Since $D(\tau)$ is reduced and $(\beta,\tau)$ is connected, $C$ cannot be contained in $\beta$, and hence $C$ also cannot be contained in $\beta'$.  Furthermore, since $C$ meets $N(\tau  + \tau' )$ in at most one point, it cannot be disjoint from $\beta\cup\beta'$.   Hence without loss of generality, we can assume that $C$ meets $\gamma$ transversely in a minimum (but non-zero) number of points.

	Suppose for the sake of contradiction that some arc of $C\cap \beta$ is disjoint from $\tau$.  Let $A$ be such an arc which is outermost in $\beta$.  Since $(\beta,\tau)$ is connected, $A$ must have both endpoints in a single component of $\gamma -  \tau$.  Now by cut and paste we can replace $A$ by an arc just outside of $\beta$ to reduce the number of points of $C\cap \gamma$.  Since this contradicts minimality, no arc of $C\cap \beta$ or $C\cap \beta'$ can be disjoint from $N(\tau  + \tau' )$.  
	 	


	
	Thus $C\cap \beta$ is a single arc $A$ meeting $\tau$ in a crossing and $C\cap \beta'=\emptyset$.  Since no arc of $C\cap\beta_i$ can intersect $N(\tau  + \tau' )$, any arc in $C\cap\beta_i$ would have both endpoints in $\gamma$, and hence could be removed by an isotopy of $C$.   It follows that the endpoints of $A$ must be contained in $\gamma- \partial(\beta_1\cup\beta_2)$.  If the endpoints of $A$ were in distinct components of $\gamma- \partial(\beta_1\cup\beta_2)$, then $C-A$ would have to intersect $\beta'\cup\beta_1\cup \beta_2$, which cannot occur.  If the endpoints of $A$ were in the same component of $\gamma- \partial(\beta_1\cup\beta_2)$, then $A$ together with an arc of $\gamma-\tau$ would be a circle meeting $D(\tau)$ transversely in one crossing, which cannot occur since $D(\tau)$ is reduced. It follows that $C$ cannot exist, and hence $N(\tau  + \tau' )$ is reduced and connected as required.

 	\medskip
 	
 	\noindent {\bf Step 2: We show that if $D(\tau)$ is alternating or positive, then $N(\tau  + \tau' )$ is alternating or positive, respectively.}

 	If $D(\tau)$ is alternating, then the last crossings before the NE and SE endpoints of $(\beta,\tau)$ are opposites (i.e., one is over and the other is under), and the last crossings before the NW and SW endpoints of $(\beta,\tau)$ are opposites.  In $N(\tau  + \tau' )$, the NE endpoint of $(\beta,\tau)$ is attached to the NW endpoint of $(\beta',\tau')$.  Also, since $(\beta',\tau')$ is a $\pi$-rotation of $(\beta,\tau)$, the last crossing before the NW endpoint of $(\beta',\tau')$ is the same as the last crossing before the SE endpoint of $(\beta,\tau)$, which is opposite the last crossing before the NE endpoint of $(\beta,\tau)$.  Using the same argument for the three other pairs of endpoints of $(\beta,\tau)$ and $(\beta',\tau')$ which are attached in $N(\tau  + \tau' )$, we see that $N(\tau  + \tau' )$ is alternating.

 	If $D(\tau)$ is positive, then we orient the strings of $(\beta,\tau)$ so that their orientation is consistent with that of $D(\tau)$, and we give the strings of $(\beta',\tau')$ the orientation obtained by rotating the oriented strands of $(\beta,\tau)$ by $\pi$.  With these orientations, all of the crossings of $(\beta,\tau)$ and hence of $(\beta',\tau')$ are positive. Now the orientations at the NE and SE endpoints of $(\beta,\tau)$ are opposites (i.e., one is inward and the other is outward), and the orientations at the NW and SW endpoints of $(\beta,\tau)$ are opposites.  Hence as in the argument for alternating, the orientations $(\beta,\tau)$ and $(\beta',\tau')$ induce a coherent orientation on $N(\tau  + \tau' )$ such that all crossings are positive.

 	\medskip

 	\noindent {\bf Step 3: We show that if $N(\tau)$ and $D(\tau)$ are both twist $A$-adequate, then $N(\tau  + \tau' )$ is twist $A$-adequate.} 
 	
  Let $S_1$ and $S_2$ denote the circles in the all-$A$ resolution of $N(\tau  + \tau' )$ which meet both $\beta$ and $\beta'$, let $\alpha_i=S_i\cap\beta$, and let $\tau_A=\alpha_1\cup \alpha_2$.

 	Now suppose for the sake of contradiction that in the all-$A$ resolution of $N(\tau  + \tau' )$, some grey arc $G$ has both endpoints on the same circle $C$.  Since $G$ comes from smoothing a crossing, it must be contained entirely in the interior of either $\beta$ or $\beta'$.  Without loss of generality, we assume that $G\subseteq \beta$, and hence $C\cap\beta$ is non-empty. Suppose that $C\subseteq \beta$.  Since the crossings of $N(\tau)$ are precisely the crossings of $N(\tau+\tau')$ that are contained in $\beta$,
 $C$ is also a circle in the all-$A$ resolution of $N(\tau)$.  This is impossible since $N(\tau)$ is $A$-adequate. Thus $C$ must be either $S_1$ or $S_2$, and hence $G$ must have both of its endpoints on either $\alpha_1$ or $\alpha_2$.  Now $\tau_A$ is contained in a circle of the all-$A$ resolution of $N(\tau)$ or $D(\tau)$.  It follows that $G$ has both of its endpoints on this circle, which is impossible since $N(\tau)$ and $D(\tau)$ are $A$-adequate.  Thus $N(\tau  + \tau' )$ is $A$-adequate as required.

 	To show that $N(\tau  + \tau' )$ satisfies the twist region condition, consider a pair of grey arcs $G_1$ and $G_2$ joining the same pair of circles $C_1$ and $C_2$ in the all-$A$ resolution of $N(\tau  + \tau' )$.  Without loss of generality, $G_1\subseteq \beta$.  Thus $C_1$ and $C_2$ both meet $\beta$.   	We consider cases according to how $C_1$ and $C_2$ intersect $\tau_A$.

	 First suppose that $C_1$ and $C_2$ are disjoint from $\tau_A$.  Since $G_1\subseteq \beta$, it follows that $C_1\cup C_2\subseteq \beta$, and hence $G_2\subseteq \beta$.   Thus
 $C_1$ and $C_2$ are also circles in the all-$A$ resolution of $N(\tau)$.  Now because $N(\tau)$ is twist $A$-adequate, $G_1$ and $G_2$ must come from the same twist region of $N(\tau)$.  Since $C_1$ and $C_2$ are disjoint from $\tau_A$, this twist region is in the interior of $\beta$, and hence must also be a twist region of $N(\tau+\tau')$.


 	Next suppose that $C_1$ is disjoint from $\tau_A$ and $C_2$ contains $\alpha_1$.   Since $G_1\subseteq \beta$, it follows that $C_1\subseteq \beta$.  Let $X$ denote the region of $\beta-\tau_A$ that contains $C_1$.  Then $G_1\cup G_2\subseteq X$, and either $\partial X$ contains $\alpha_1\cup \alpha_2$ and is disjoint from $\partial \beta_1\cup\partial \beta_2$, or $\partial X$ is disjoint from $\alpha_2$ and from the top and bottom arcs of $\partial \beta-\partial \tau_A$.  First we suppose the former.  Let $C_3$ denote the circle of  the all-$A$ resolution of $D(\tau)$ that contains  $\alpha_1$.  Then $G_1$ and $G_2$ go between $C_1$ and $C_3$.   Since $D(\tau)$ is twist $A$-adequate, this means that $G_1$ and $G_2$ come from crossings in the same twist region of $D(\tau)$, and the twist region is contained in $X$.  Now since $X$ is disjoint from $\partial \beta_1\cup\partial \beta_2$, this twist region of $D(\tau)$ is in the interior of $\beta$, and hence must also be a twist region of $N(\tau+\tau')$.
	

Next suppose that $\partial X$ is disjoint from $\alpha_2$ and from the top and bottom arcs of $\partial \beta-\partial \tau_A$.  Let $C_3$ denote the circle of  the all-$A$ resolution of $N(\tau)$ containing $\alpha_1\cup \alpha_2$.  Then $G_1$ and $G_2$ go between $C_1$ and $C_3$.   Since $N(\tau)$ is twist $A$-adequate, this means that $G_1$ and $G_2$ come from crossings in the same twist region of $N(\tau)$, and the twist region is contained in $X$.  Now since $X$ is disjoint from the top and bottom arcs of $\partial \beta-\partial \tau_A$, this twist region of $N(\tau)$ is in the interior of $\beta$, and hence must also be a twist region of $N(\tau+\tau')$.

 	Finally, suppose that $\alpha_1\subseteq C_1$ and $\alpha_2\subseteq C_2$.  Since $G_1\subseteq \beta$, this means that $G_1$ goes between $\alpha_1$ and $\alpha_2$.   Thus in the all-$A$ resolution of $N(\tau)$, $G_1$ has both of its endpoints on the circle consisting of $\tau_A$ and the top and bottom arc of $\partial \beta-\partial \tau_A$.  But this cannot occur since $N(\tau)$ is $A$-adequate.  It follows that $N(\tau  + \tau' )$ is twist $A$-adequate, and hence we have proven the lemma.  \end{proof}

 \begin{lemma}  \label{knottedarc} Let $(B,T)$ be a tangle whose projection $(\beta,\tau)$ is not prime and contains at least one crossing.  If $D(\tau)$ is reduced, and either alternating, positive, or semi-adequate.  Then there is a disk $\delta\subseteq\beta$ whose intersection with $\tau$ is the projection of a knotted arc in $T$.
 \end{lemma}

 \begin{proof}
 Recall that since $(\beta,\tau)$ is not prime, by definition there is a circle $C\subseteq \beta$ meeting $\tau$ in two points such that $C$ bounds a disk $\delta\subseteq\beta$ whose intersection with $\tau$ is an arc with at least one crossing.  We need to know that the arc $\tau\cap \delta$ is actually the projection of a knotted arc.
 
 Let $S$ be a sphere obtained by capping off $C$ in $B$ above and below the arcs of $T$ such that the ball $\Delta$ bounded by $S$ in $B$ projects to $\delta$ and the arc $T\cap \Delta$ projects to the arc $\tau\cap \delta$.  Let $K$ be obtained from the arc $T\cap \Delta$ by joining its endpoints in $\partial \Delta$ such that all of the crossings in its projection $\kappa$ are in $\tau\cap \delta$.  Since $D(\tau)$ is reduced, and either alternating, positive, or semi-adequate, then so is $\kappa$.    It follows that $K$ is a non-trivial knot when $\kappa$ is alternating by Kauffman, Murasugi, and Thistlethwaite's Theorem, when $\kappa$ is positive by Stoimenow's Theorem, and when $\kappa$ is semi-adequate by Thistlethwaite's Theorem.  Thus $T\cap \Delta$ is a knotted arc, and hence $\tau\cap \delta$ is the projection of a knotted arc.
 \end{proof}

  Next we bring together the hypotheses of Menasco {\bf (M)}, Ozawa {\bf(O)}, and Futer et al. {\bf(F)} in the following definition.

 \begin{defi}  Let $(\beta,\tau)$ be a projection of a $2$-string tangle.  If both $N(\tau)$ and $D(\tau)$ are reduced and at least one of the following holds, then we say that $(\beta,\tau)$ satisfies {\bf MOF}. 
 	\medskip
 	
 	\noindent {\bf(M)}   At least one of $N(\tau)$ or $D(\tau)$ is alternating.  
 	
 	\medskip
 	
 	\noindent{\bf(O)}  At least one of $N(\tau)$ or $D(\tau)$ is positive.  
 	
 	\medskip
 	
 	\noindent{\bf(F)}  $N(\tau)$ and $D(\tau)$ are connected and both are twist $A$-adequate or both are twist $B$-adequate. 
 	
 \end{defi}


We now use the above lemmas to prove the following.

 \begin{theorem}\label{tangle} Let $(B,T)$ be a non-rational tangle whose projection $(\beta,\tau)$ satisfies MOF.  Then $(\beta,\tau)$ is prime if and only if $(B, T)$ is prime.\end{theorem}
 
 Note that if $(\beta,\tau)$ is strongly alternating and connected, it follows from Lickorish and Thistlethwaite \cite{LT} that $(B,T)$ cannot be rational.

 \begin{proof} Suppose that $(\beta,\tau)$ is not prime.  If hypothesis (M) or (O) is satisfied, we can replace $(\beta,\tau)$ by a $\frac{\pi}{2}$ rotation if necessary so that $D(\tau)$ is alternating or positive.  It follows that $(\beta,\tau)$ satisfies the hypotheses of Lemma~\ref{knottedarc}.  Thus there is a disk $\delta\subseteq\beta$ whose intersection with $\tau$ is the projection of a knotted arc in $T$.  It follows that $(B,T)$ is not prime.

 	To prove the converse, suppose that $(B,T)$ is not prime.  If $(\beta,\tau)$ is disconnected, then there is an arc $A\subseteq \beta-\tau$ which separates the strings of $\tau$.  Let $A_1$ and $A_2$ be the arcs of $\partial \beta-A$.  Then for each $i$, let $C_i=A\cup A_i$ and $\beta_i$ be the disk in $\beta$ bounded by $C_i$.  Then one string of $\tau$ is in each $\beta_i$.  If both of these strings are trivial, then $(\beta,\tau)$ would be trivial and hence would be rational.  Thus at least one of these strings is non-trivial, and so $(\beta,\tau)$ is not prime.  
	
	Thus we assume that $(\beta,\tau)$ is connected.  Let $(\beta', \tau')$ be a $\pi$-rotation of $(\beta,\tau)$, let $(B',T')$ be a tangle whose projection is $(\beta', \tau')$, let $L=N(T+T')$ with projection $G=N(\tau+ \tau')$.   Since $(\beta,\tau)$ is connected and satisfies MOF,  Lemma~\ref{lem:cp} implies that $G$ is reduced and connected, and is either alternating, positive, or twist-adequate.  Furthermore, if $G$ is alternating, then by the first part of Menasco's Theorem, $L$ is non-split.
	
	Since $(B,T)$ is neither prime nor rational, it must be composite. Thus there is a ball in $B$ meeting $T$ in a knotted arc.  As the same is true for $(B',T')$, the link $L$ is composite.   Now we can apply the theorem of Menasco, Ozawa, or Futer et al. to conclude that $G$ is not prime.   Hence there is a circle $C$ in the sphere of projection $P$ meeting $G$ transversely in two points such that neither component of $P-C$ meets $G$ in a trivial arc.   If $C$ were contained in $\beta$ or $\beta'$, then $(\beta,\tau)$ would not be prime and we would be done.   
	
	We suppose that is not the case.  Since $C$ meets $G$ in just two points, $C$ cannot be contained in $P-(\beta\cup \beta')$.  Thus without loss of generality, we assume that $C$ meets $\partial \beta$ in a minimal but non-zero number of points.   Let $A$ be an arc of $C\cap \beta$ which is outermost in $\beta$.  If $A$ were disjoint from $\tau$ then, since $(\beta,\tau)$ is connected, $A$ could be isotoped into $\partial \beta$, and hence $C\cap \partial \beta$ would not be minimal.  Thus every arc of $C\cap \beta$ must intersect $\tau$.

 	Suppose that one component of $\beta-A$ is a disk $D$ that meets $\tau$ in a single arc.  If $D\cap \tau$ were a non-trivial arc then $(\beta,\tau)$ would not be prime and again we would be done.  Thus we assume that one component of $\beta-A$ is a disk $D$ that meets $\tau$ in a single trivial arc.   If $A$ intersects that string of $\tau$ in two points, then $C-A$ must be disjoint from $G$.  Hence the component of $P-C$ containing $D$ would intersect $G$ in a trivial arc, contrary to the definition of $C$. But if $A$ intersects $\tau$ in just one point, then we can assume $A$ is outermost on $\beta$ with respect to this property.  Now since $D\cap \tau$ is a trivial arc, we can replace $A$ by an arc just outside of $\beta$ also meeting $G$ in just one point.  As this contradicts the minimality of $C\cap \partial \beta$, we can assume that each component of $\beta-A$ meets $\tau$ in two arcs.

It follows that each component of $\partial \beta-\partial A$ contains two endpoints of $T$, which may or may not be endpoints of the same string.  Now (referring to the strips $\beta_1$ and $\beta_2$ from Figure 6) either one endpoint of $A$ is in $\partial\beta_1$ and the other endpoint is in $\partial\beta_2$ or the two endpoints of $A$ are in distinct components of $\partial \beta-(\partial \beta_1\cup \partial\beta_2)$.  But since $A$ meets $G$ in two points, $C-A$ must be disjoint from $G$, which is impossible since $(\beta,\tau)$ and $(\beta',\tau')$ are connected.  It follows that $C$ is contained in $\beta$ or $\beta'$, and hence $(\beta,\tau)$ is not prime.\end{proof}

Now we prove Cromwell's Theorem as a corollary to the above theorem and its proof.

\begin{cor}   Let $(B,T)$ be a tangle whose projection $(\beta,\tau)$ is strongly alternating. Then $(B, T)$ is prime if and only if $(\beta,\tau)$ is prime.
\end{cor}

\begin{proof}  Since $(\beta,\tau)$ is strongly alternating, $(\beta,\tau)$ satisfies MOF.  Observe that neither Lemma~\ref{knottedarc} nor the first two paragraphs of the proof of Theorem~\ref{tangle} require that $(B, T)$ be non-rational.  It follows that if $(\beta,\tau)$ is not prime then $(B,T)$ is not prime, and if $(B,T)$ is not prime and $(\beta,\tau)$ is disconnected then $(\beta,\tau)$ is not prime.

Thus we assume that $(B,T)$ is not prime and $(\beta,\tau)$ is connected.  Now by Lickorish and Thistlethwaite \cite{LT}, since $(\beta,\tau)$ is strongly alternating and connected, $(B,T)$ cannot be rational.  Hence we can apply Theorem~\ref{tangle} to conclude that $(\beta,\tau)$ is not prime.
\end{proof}

 

 \section{Splittings of spatial graphs}\label{Graphs}

By a {\bf spatial graph} we mean a graph embedded in $S^3$, where we allow graphs to have edges which are loops and multiple edges going between the same pair of vertices.  If an edge of a graph is a loop, then the two pieces of the loop adjacent to the vertex are treated as separate edges.  

\begin{figure}[h!]
			\centering	\includegraphics[width=.65\textwidth]{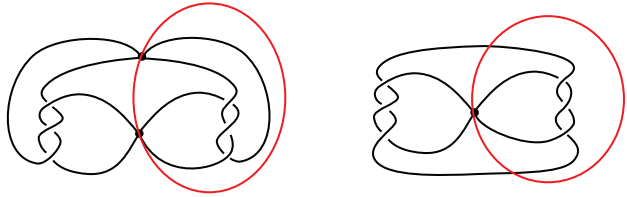}
			\caption{Sawollek considers splitting spheres and circles which pass through vertices. }
			\label{spheres}
			\end{figure}

Menasco's Theorem has been extended to reduced alternating spatial graphs independently by Sawollek \cite{Sawollek2} and Adams et al. \cite{Adams} (using different definitions).  However, Sawollek considers splitting spheres and circles which pass through vertices, while Adams et al. does not (see Figure~\ref{spheres}).

 Below are Sawollek's definitions.   

\begin{defi}\label{spatialSplit}  Let $\Gamma$ be a $4$-valent spatial graph. 
\begin{enumerate}
\item $\Gamma$ is {\bf split} if there is a sphere $\subseteq S^3-\Gamma$ which separates $\Gamma$.   

\item Suppose that $\Gamma$ is not split.  Then $\Gamma$ is {\bf vertex split} if there is a sphere $S\subseteq S^3$ intersecting $\Gamma$ in a non-zero number of vertices such that each component of $S^3-S$ contains two  edges incident to each vertex in $S\cap \Gamma$.    

\item Suppose that $\Gamma$ is neither split nor vertex-split.  Then $\Gamma$ is {\bf composite} if there is a sphere $S\subseteq S^3$ which meets $\Gamma$ transversely in any number (possibly zero) of vertices and exactly two points in the interiors of edges such that neither component of $S^3-S$ meets $\Gamma$ can be isotoped into $S$, keeping $S$ pointwise fixed.  Otherwise, $\Gamma$ is {\bf prime}.  \end{enumerate}\end{defi}

\begin{defi}\label{projectionSplit} Let $G$ be a projection of a $4$-valent spatial graph on a sphere of projection $P$.  

\begin{enumerate}
\item  $G$ is {\bf split} if there is a circle $C\subseteq P-G$ which separates $G$.  

\item Suppose that $G$ is not split.  Then $G$ is {\bf vertex split} if there is a circle $C\subseteq P$ which intersects $G$ transversely in some non-zero number of vertices (i.e., such that each component of $P-C$ contains two of the edges incident to each of the vertices in $C\cap G$).  

\item Suppose that $G$ is neither split nor vertex-split.  Then $G$ is {\bf composite} if there is a circle $C\subseteq P$ which meets $G$ transversely in any number (possibly zero) of vertices and exactly two points in the interiors of edges such that neither component of $P-C$ meets $G$ in subgraphs with no crossings. Otherwise $G$ is {\bf prime}. 

\end{enumerate}\end{defi}

In fact, we have corrected Sawollek's definitions of composite for spatial graphs and their projections, since the original definitions do not require that the components of the split apart graph or of its projection be non-trivial.  In this case, every $4$-valent spatial graph and its projection which are not split or vertex split are composite because you can always choose the splitting sphere or splitting circle to cut off a small trivial arc in a single edge.  This is similar to the issue we mentioned with Menasco's original statement of his theorem.

Below we give Sawollek's definition of reduced and alternating and his generalization of  Menasco's Theorem to $4$-valent spatial graphs. 

\begin{defi} \label{RedAltSawollek} Let $G$ be a projection of a $4$-valent spatial graph.    We say that $G$ is {\bf reduced} and {\bf alternating}, if every smoothing of the vertices in the sphere of projection yields a reduced alternating projection of a link.   \end{defi}

 \begin{wrongM}\cite{Sawollek2} \label{SawollekM} Let $G$ be a reduced alternating projection of a $4$-valent spatial graph $\Gamma$.  Then 
  \begin{enumerate}
 \item $\Gamma$ is split if and only if $G$ is split.
 
 \item $\Gamma$ is vertex split if and only if $G$ is vertex split.
 
\item $\Gamma$ is composite if and only if $G$ is composite.
\end{enumerate}\label{wrongM}
 \end{wrongM}

\begin{figure}[h!]
			\centering	\includegraphics[width=.8\textwidth]{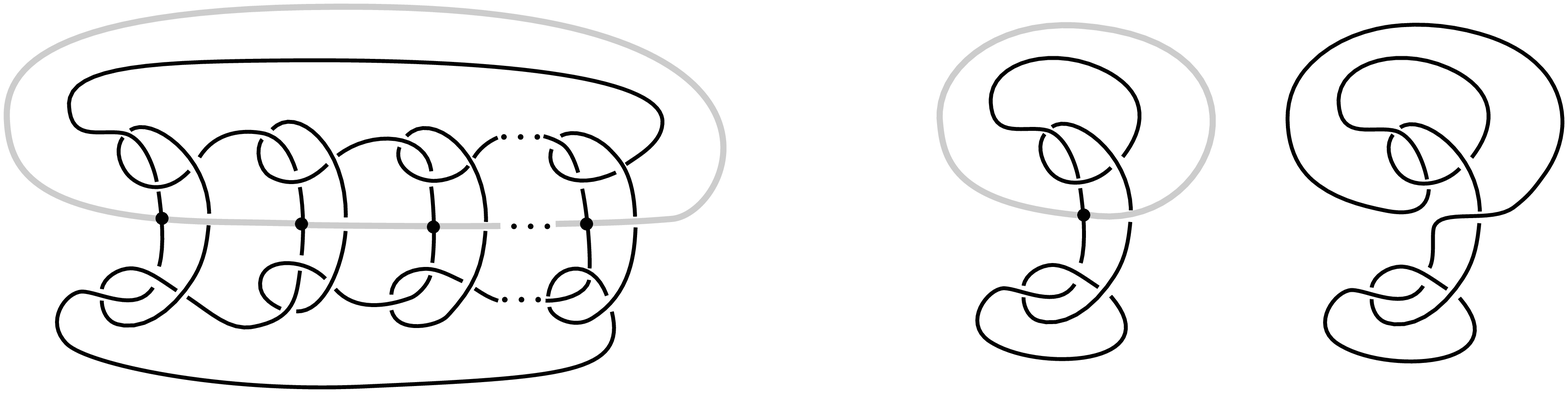}
			\caption{When there is more than one vertex, this spatial graph is vertex split but its reduced alternating projection is not. If there is only one vertex, the projection is not reduced. }
			\label{fig:Erica}
			\end{figure}

\begin{figure}[h!]
			\centering	\includegraphics[width=.6\textwidth]{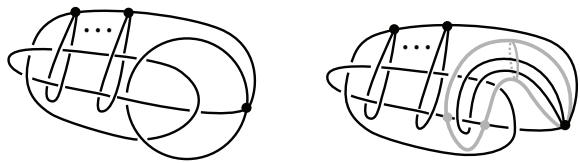}
			\caption{This spatial graph is composite, but its reduced alternating projection is not.}
			\label{notprime}
			\end{figure}

Parts (2) and (3) are noteworthy because they concern spheres which pass through vertices, and hence don't appear to follow closely from results or techniques for knots and links.  However, although Sawollek's proof of Part (1) is correct and the backwards direction of all three parts are easily seen to be true, the result fails for the forward direction of Part (2) for $4$-valent graphs with more than one vertex and for the forward direction of Part (3) when the splitting sphere contains at least one vertex.   The problem is that the proof assumes that when a splitting sphere $S$ meets a vertex, the two edges which are on the same side of $S$ are adjacent in the projection.  This is not true in general, as illustrated in Figures~\ref{fig:Erica} and \ref{notprime}.  


 In particular, for the projection of the $4$-valent graph with $n>1$ vertices in Figure~\ref{fig:Erica}, there is a vertex splitting sphere passing through all of the vertices bounding a ball that contains precisely the two grey edges at each vertex.  However, this projection is reduced and alternating, and has no vertex splitting circle.   When there is only one vertex, the graph does not satisfy the hypotheses of Sawollek's Theorem since, as shown on the right, one of the smoothings of the vertex yields a link is not reduced.

Figure~\ref{notprime} shows that, regardless of the number of vertices, the forward direction of Part (3) of Sawollek's Theorem is false if the splitting sphere contains a vertex.  In particular, the projection on the left is reduced, alternating, and prime.  However, as illustrated on the right, there is a grey sphere $S$ meeting the graph in one vertex together with two points in the interior of an edge. The sphere $S$ splits the graph into two pieces such that neither component of $S^3-S$ meets $\Gamma$ in one or two arcs which can be isotoped into $S$.  Thus the spatial graph is composite, though its reduced alternating projection is not.

 In Section~\ref{OurGraphs}, we prove a stronger version of the forward direction of Part (2) in our Theorem 2, and of Part (3) in our Theorem 3 for graphs with a single $4$-valent vertex, when the splitting sphere is disjoint from the vertex.

\medskip

 Next we consider Adams et al.'s generalization of Menasco's Theorem, which is proved using the machinery of Menasco.   Adams et al.'s results are distinct from those of Sawollek not only because the definitions are different and they allow vertices of any valence, but because they do not allow spheres to cut through vertices.  The following definition is given in Adams et al.  \cite{Adams}.			
			
\begin{defi}\label{RedAltAdams}  Let $G$ be a projection of a spatial graph.  
\begin{enumerate}
\item $G$ is {\bf reduced} if there is no circle $C$ in the sphere of projection such that $C$ intersects $G$ transversely in a single crossing. 
\item A {\bf region} $R$ of $G$ is a connected component of $S^2-G$ with $\partial R\subseteq G$.  
\item Given a region $R$ of $G$, a {\bf segment} $E$ of $R$ is a connected subset of $\partial R$ whose interior intersects no crossings of $G$, such that $\partial E$ is a non-empty set of crossings.  
\item An {\bf alternating} segment is a segment which has precisely two endpoints, where one endpoint is an undercrossing and the other endpoint is an overcrossing. 
\item A {\bf alternating region} is a region such that every segment of it is alternating. \end{enumerate} \end{defi}

\begin{figure}[h!]
			\centering	\includegraphics[width=0.65\textwidth]{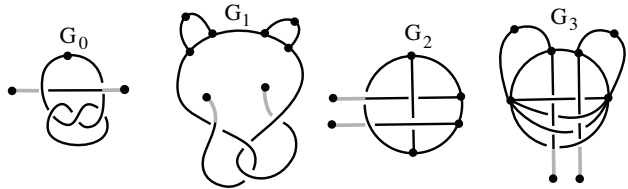}
			\caption{The grey arcs are not segments, hence these graph projections are alternating by the above definition.   }
			\label{valence1}
			\end{figure}

According to this definition, no piece of a graph containing a valence $1$ vertex can be a segment.  For example, the grey arcs  in Figure~\ref{valence1} are not segments.  Thus all of the segments in  Figure~\ref{valence1} are alternating, and hence all of the regions are as well.   Adams et al. also defines an alternating projection and states that a projection is alternating if and only if all of its regions are alternating. It follows that the projections in Figure~\ref{valence1} are reduced and alternating.
  
However, Adams has informed us that the definitions of a region and a segment given in \cite{Adams} were not what the authors had intended.  Rather, they had intended to define a region $R$ as a connected component of the complement of a {\it regular neighborhood} of $G$, and a segment of $R$ as an interval in the boundary of a {\it regular neighborhood} of $G$ containing no crossings and bounded by crossings on its ends.  In Figure~\ref{valence1}, the boundary of a regular neighborhood of each of the grey arcs is a U-shaped arc around the vertex of valence $1$ whose endpoints are a pair of identical crossings.  With this definition of a segment the grey arcs in Figure~\ref{valence1} are segments which are not alternating, and hence the projections are not alternating.  In fact, with this definition, no connected spatial graph projection with a valence $1$ vertex and at least one crossing can be alternating.

\begin{defi}\label{Adams} Let $\Gamma$ be a spatial graph, let $F$ be a sphere which meets $\Gamma$ transversely in $n$ points in the interiors of edges, and let $B$ be a component of $S^3-F$.  If $B\cap \Gamma$ contains either a vertex of valence at least $3$ or a subgraph which cannot be isotoped into $F$ while fixing its boundary points, then we say $B\cap \Gamma$ is {\bf non-trivial}.  If both components of $(S^3-F)\cap \Gamma$ are non-trivial and $n$ is minimal, then we say $\Gamma$ is {\bf $n$-composite}.  \end{defi}

 \begin{AdamsM}\label{AdamsMenasco}{\bf (Theorem 3.1 of \cite{Adams})}  Let $G$ be a reduced alternating projection of an $n$-composite spatial graph with $n<4$. Then there is a circle $C$ meeting $G$ transversely in $n$ points such that each component of $S^2-C$ either contains a vertex of valence at least $3$ or at least one crossing.    \end{AdamsM}

Note that this is a corrected version of the statement in \cite{Adams}, which had a similar issue to the one that was in Menasco's original statement.  

 By contrast with this Theorem, each spatial graph in Figure~\ref{valence1} is $n$-composite, yet according to the definitions in \cite{Adams} the $G_n$ are reduced and alternating and there is no circle $C$ meeting $G_n$ in $n$ points such that each component of $S^2-C$ either contains a vertex of valence at least $3$ or at least one crossing.  On the other hand, using the definitions that Adams et al. had intended, the projections in Figure~\ref{valence1} would not be alternating and hence the above theorem would not apply.

\medskip

\begin{figure}[h!]
			\centering	\includegraphics[width=0.65\textwidth]{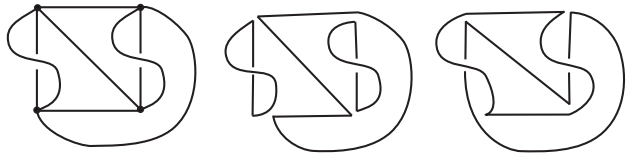}
			\caption{The projection on the left is reduced and alternating according to Adams et al., but is neither reduced nor alternating according to Sawollek. }
			\label{RedAlt}
			\end{figure}

We close this section by comparing some of the definitions for $4$-valent graphs in Sawollek with those in Adams et al.. In particular, Sawollek's definitions of reduced and alternating are stronger than than those of Adams et al., since if all smoothings of the vertices yield a reduced alternating projection, then there cannot be a circle meeting the graph projection in a single crossing nor a region of the graph projection which is not alternating.  For example, the projection on the left in Figure~\ref{RedAlt} is reduced and alternating according to Adams' et al., but the smoothings of vertices in the center and right show that it is neither reduced nor alternating according to Sawollek.  

\begin{figure}[h!]
			\centering	\includegraphics[width=.4\textwidth]{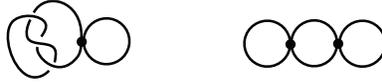}
			\caption{These spatial graphs are vertex split hence not composite according to Sawollek, but are $2$-composite according to Adams et al.. }
			\label{CompositeEight}
			\end{figure}

Neither Sawollek's definition of composite nor Adams et al.'s definition of $2$-composite implies the other.  For example, the graphs in Figure~\ref{CompositeEight} are $2$-composite according to Adams et al., but are vertex split and hence not composite according to Sawollek.  In addition, for any sphere which meets one of these spatial graphs in two non-vertex points, on one side of the sphere the subgraph can be isotoped into the sphere.  Thus the sphere does not split the graph into two non-trivial subgraphs in the sense of Sawollek, though such a sphere does split it into two non-trivial graphs in the sense of Adams et al..

On the other hand, the spatial graph in Figure~\ref{notprime} is composite according to Sawollek, but a splitting sphere passes through the vertex, and hence it is not $2$-composite according to Adams et al..

\section{Our results for figure eight graphs}\label{OurGraphs}
         
In this section, we consider figure eight graphs (i.e.,  $4$-valent graphs with a single vertex), and prove the forward direction of Part (2) of Sawollek's Theorem and the forward direction of Part (3) when the splitting sphere is disjoint from the vertex.  Recall that Figure~\ref{fig:Erica} shows that the forward direction of Part (2) is false for $4$-valent graphs with more than one vertex and Figure~\ref{notprime} shows that the forward direction of Part (3) is false if the splitting sphere contains vertices.

\begin{figure}[h!]
			\centering	\includegraphics[width=.4\textwidth]{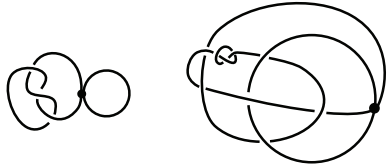}
			\caption{These embeddings contain local knots.  The one on the left is vertex split and the one on the right is composite. }
			\label{LocalKnot}
			\end{figure}

 We will say that a spatial graph $\Gamma$ contains a {\bf local knot} if there is a sphere $S$ meeting $\Gamma$ in two points in the interior of an edge such that one component of $(S^3-S)\cap \Gamma$ is a knotted arc.  The embeddings in Figure~\ref{LocalKnot} each contain a local knot and have a reduced alternating projection.  However, the embedding and projection on the left are vertex split, and hence not compositive.  If $\Gamma$ contains a local knot and is not vertex split, then $\Gamma$ is composite.

\begin{theorem}\label{alternating}  Let $\Gamma$ be an embedding of a vertex split figure eight graph whose projection $G$ is reduced, alternating, and contains at least one crossing.  Then there is a local knot in $\Gamma$ whose projection is contained in a disk with boundary meeting the vertex transversely.


\end{theorem}

 Note that this theorem is stronger than Part (2) of Sawollek's Theorem for figure eight graphs, since it tells us that a reduced alternating vertex split graph contains a local knot which you can see in the projection.  It is also independent of  Adams et al.'s results because Adams et al. does not consider spheres which pass through vertices.

\begin{proof}  Let $N(v)\subseteq S^3$ be a neighborhood of $v$ whose projection $n(v)$ contains no crossings of $G$.  Let $B=\mathrm{cl}(S^3-N(v))$, $T=B\cap \Gamma$, and $(\beta,\tau)$ be the projection of the tangle $(B,T)$.  Since $G$ is reduced, alternating, and contains at least one crossing, $(\beta,\tau)$ is strongly alternating.

Suppose that $(\beta,\tau)$ is disconnected. Then there is an arc $\alpha$ which separates the strings of $\tau$.  Hence without loss of generality, $N(\tau)$ is a disconnected projection of the link $N(T)$.  Since $G$ contains at least one crossing, the projection $K$ of one of the components of $N(T)$ contains at least one crossing.  Furthermore, because $N(\tau)$ is reduced and alternating and there are no crossings between the strings of $\tau$, $K$ is reduced and alternating.  Thus by Kauffman, Murasugi, and Thistlethwaite's Theorem, $K$ is the projection of a non-trivial knot.   Now we extend the arc $\alpha$ into $n(v)$ by adding an arc through $v$ to obtain a circle $C$ in the sphere of projection which separates the loops of $G$ such that $C$ bounds a disk containing the projection of a non-trivial knot.  

Next we assume for the sake of contradiction that $(\beta,\tau)$ is connected.  Let $S$ be a vertex splitting sphere for $\Gamma$.  Then $S$ separates the loops of $\Gamma$.  Thus, by making $N(v)$ smaller if necessary, we can assume that $S$ meets $N(v)$ in a disk intersecting $\Gamma$ only in $v$ and $S$ meets $B$ in a disk $\Delta$ which separates the strings of $(B,T)$.


By Lickorish and Thistlethwaite \cite{LT}, since $(\beta,\tau)$ is strongly alternating and connected, $(B,T)$ cannot be rational. Thus $(B,T)$ must be composite, and hence we can apply Theorem~\ref{tangle} to obtain a circle $C$ bounding a disk $D_1\subseteq \beta$ whose intersection with $\tau$ is a non-trivial arc $\alpha$.  Now there is a ball $B_1\subseteq B$ projecting to $D_1$ such that $B_1\cap T$ is an arc $A$ which projects to $\alpha$.  Let $A_1$ be an arc in $\partial B_1$ with the same endpoints as $A$ such that $A_1$ injectively projects to an arc $\alpha_1$ of $C$, and let $(B,T_1)$ be obtained from $(B,T)$ by replacing $A$ by $A_1$.  Then the projection $(\beta, \tau_1)$ of $(B,T_1)$ is obtained from $(\beta,\tau)$ by replacing $\alpha$ by $\alpha_1$, and hence $(\beta, \tau_1)$ has fewer crossings in one of the strings than $(\beta,\tau)$ had.   

Observe that since the splitting disk $\Delta$ is disjoint from the strings of $T$ and $B_1$ meets $T$ in an arc of a single string, any circles of intersection of $\Delta$ and $\partial B_1$ can be removed to obtain an isotopic splitting disk $\Delta_1$ for $(B,T)$ which is disjoint from $B_1$.  It follows that $\Delta_1$ is also a splitting disk for $(B,T_1)$.  Since $\alpha$ is a non-trivial arc of $N(\tau)$ and $D(\tau)$ which are alternating, its last crossing adjacent to one of its ends must be an overcrossing and its last crossing adjacent to its other end must be an undercrossing.  Hence $N(\tau_1)$ and $D(\tau_1)$ must also be alternating.  Furthermore, any nugatory crossing of $\tau_1$ would be disjoint from $D_1$ and hence would also be a nugatory crossing of $\tau$. Thus $(\beta, \tau_1)$ is strongly alternating.   Finally, since the strings of $\tau_1$ cross in the same points as those of $\tau$ did, $(\beta, \tau_1)$ is connected.

Thus we can repeat the above two paragraphs with $(B,T_1)$ in place of $(B,T)$ to obtain a new tangle $(B,T_2)$ whose projection $(\beta, \tau_2)$ is connected, strongly alternating, and has fewer crossings in one of the strings than $(\beta, \tau_1)$ had.  We repeat this process, reducing the number of crossings within one of the strings at each stage.  But we cannot continue this indefinitely since there are only finitely many crossings within each string.  Thus in fact, $(\beta,\tau)$ could not be connected.
\end{proof}


\begin{theorem}\label{figureeight}  Let $\Gamma$ be an embedding of a figure eight graph with vertex $v$ and projection $G$ such that $\Gamma$ contains a local knot.  Let $N(v)$ be a neighborhood of $v$ whose projection $n(v)$ contains no crossings of $G$, $B=\mathrm{cl}(S^3-N(v))$, and $T=B\cap \Gamma$, and suppose the projection $(\beta,\tau)$ of $(B,T)$ satisfies MOF.  Then there is a disk in $ \beta$ whose intersection with $G$ is the projection of a knotted arc.  
\end{theorem}

\begin{figure}[h!]
			\centering	\includegraphics[width=.5\textwidth]{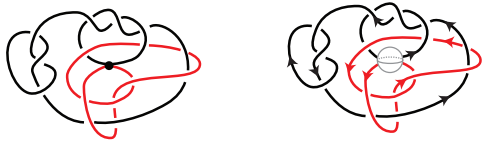}
			\caption{We can apply Theorem ~\ref{figureeight} to this embedded graph, but Sawollek's Theorem and Adams et al.'s Theorem do not apply. }
			\label{Positive}
			\end{figure}

Note that for figure eight graphs, Theorem~\ref{figureeight} is stronger than Part (3) of Sawollek's Theorem and Adams et al.'s Theorem because it applies to embeddings that are not necessarily alternating.  For example, Theorem~\ref{figureeight} applies to the embedded graph in Figure~\ref{Positive}, where we have removed a neighborhood of the vertex to obtain a $2$-string tangle $(B,T)$ with projection $(\beta,\tau)$.  Observe that $N(\tau)$ and $D(\tau)$ are reduced and positive, but not alternating.  

\begin{proof}   Since $\Gamma$ contains a local knot, there is a sphere $S$ meeting $\Gamma$ in two points in the interior of an edge such that one component of $S^3-S$ meets $\Gamma$ in a knotted arc. By making $N(v)$ smaller if necessary, we can assume that $S$ is disjoint from $N(v)$.  It follows that the component of $S^3-S$ which meets $\Gamma$ in a knotted arc is contained in $B$, and hence $(B,T)$ is not prime.  

Now suppose that $(B,T)$ is rational.  Then the instructions for how to move the endpoints of $T$ on $\partial B$ to obtain a trivial tangle, give us instructions for how to move the endpoints of the edges of $\Gamma$ around $v$ to make $\Gamma$ planar. But this is impossible since $\Gamma$ contains a local knot.  Hence $(B,T)$ cannot be rational.  

Thus since $(\beta,\tau)$ satisfies MOF, we can apply Theorem~\ref{tangle} to conclude that $(\beta,\tau)$ is not prime.  Finally, by Lemma~\ref{knottedarc}, there is a disk in $\beta$ whose intersection with $G$ is the projection of a knotted arc.  \end{proof}

It remains an open question of whether Theorem~\ref{alternating} can be extended to $4$-valent graphs which have more than one vertex.

\section*{Acknowledgement}The authors wish to thank Colin Adams for reading a preliminary version of this paper and for helpful discussions.


\bibliographystyle{amsplain}

  \bibliography{SplittingArxivSubmitted.bib}

\providecommand{\bysame}{\leavevmode\hbox to3em{\hrulefill}\thinspace}
\providecommand{\MR}{\relax\ifhmode\unskip\space\fi MR }
\providecommand{\MRhref}[2]{%
  \href{http://www.ams.org/mathscinet-getitem?mr=#1}{#2}
}
\providecommand{\href}[2]{#2}
\begin{thebibliography}{10}

\bibitem{Adams}
Colin Adams, Ryan Dorman, Kerryann Foley, Jonathan Kravis, and Sam Payne,
  \emph{Alternating graphs}, J. Combin. Theory Ser. B \textbf{77} (1999),
  no.~1, 96--120. \MR{1710534}

\bibitem{CromHomo}
P.~R. Cromwell, \emph{Homogeneous links}, J. London Math. Soc. (2) \textbf{39}
  (1989), no.~3, 535--552. \MR{1002465}

\bibitem{Crom}
\bysame, \emph{Lonely knots and tangles: identifying knots with no companions},
  Mathematika \textbf{38} (1991), no.~2, 334--347 (1992). \MR{1147833}

\bibitem{Futer}
David Futer, Efstratia Kalfagianni, and Jessica~S. Purcell, \emph{Hyperbolic
  semi-adequate links}, Comm. Anal. Geom. \textbf{23} (2015), no.~5, 993--1030.
  \MR{3458811}

\bibitem{Kauff}
Louis~H. Kauffman, \emph{State models and the {J}ones polynomial}, Topology
  \textbf{26} (1987), no.~3, 395--407. \MR{899057}

\bibitem{LT}
W.~B.~R. Lickorish and M.~B. Thistlethwaite, \emph{Some links with nontrivial
  polynomials and their crossing-numbers}, Comment. Math. Helv. \textbf{63}
  (1988), no.~4, 527--539. \MR{966948}

\bibitem{Menasco}
W.~Menasco, \emph{Closed incompressible surfaces in alternating knot and link
  complements}, Topology \textbf{23} (1984), no.~1, 37--44. \MR{721450}

\bibitem{Mura}
Kunio Murasugi, \emph{Jones polynomials and classical conjectures in knot
  theory}, Topology \textbf{26} (1987), no.~2, 187--194. \MR{895570}

\bibitem{Ozawa3}
Makoto Ozawa, \emph{Closed incompressible surfaces in the complements of
  positive knots}, Comment. Math. Helv. \textbf{77} (2002), no.~2, 235--243.
  \MR{1915040}

\bibitem{Sawollek2}
J\"{o}rg Sawollek, \emph{Alternating diagrams of {$4$}-regular graphs in
  {$3$}-space}, Topology Appl. \textbf{93} (1999), no.~3, 261--273.
  \MR{1688477}

\bibitem{Shimo}
Koya Shimokawa, \emph{Parallelism of two strings in alternating tangles}, J.
  Knot Theory Ramifications \textbf{7} (1998), no.~4, 489--502. \MR{1633011}

\bibitem{Stoi}
Alexander Stoimenow, \emph{Positive knots, closed braids and the {J}ones
  polynomial}, Ann. Sc. Norm. Super. Pisa Cl. Sci. (5) \textbf{2} (2003),
  no.~2, 237--285. \MR{2004964}

\bibitem{Thistle}
Morwen~B. Thistlethwaite, \emph{A spanning tree expansion of the {J}ones
  polynomial}, Topology \textbf{26} (1987), no.~3, 297--309. \MR{899051}

\bibitem{ThistlethwaiteAdequate}
\bysame, \emph{On the {K}auffman polynomial of an adequate link}, Invent. Math.
  \textbf{93} (1988), no.~2, 285--296. \MR{948102}

\bibitem{This}
\bysame, \emph{On the {K}auffman polynomial of an adequate link}, Invent. Math.
  \textbf{93} (1988), no.~2, 285--296. \MR{948102}

\bibitem{Tsukamoto}
Tatsuya Tsukamoto, \emph{A criterion for almost alternating links to be
  non-splittable}, Math. Proc. Cambridge Philos. Soc. \textbf{137} (2004),
  109--133.

\end{thebibliography}

\end{document}